\documentclass[final]{elsarticle}  

\usepackage{amssymb,latexsym}
\usepackage{amsmath}

\biboptions{sort&compress}

\journal{: \; Int. T. Special Functions (ACCEPTED)}

\begin{document}

\newtheorem{teo}{Theorem}
\newproof{prova}{Proof}

\begin{frontmatter}

\title{An Euler-type formula for $\beta(2n)$ and closed-form expressions for a class of zeta series}

\author{F. M. S. Lima}

\address{Institute of Physics, University of Brasilia, P.O. Box 04455, 70919-970, Brasilia-DF, Brazil}


\ead{fabio@fis.unb.br}


\begin{abstract}
In a recent work, Dancs and He found an Euler-type formula for $\,\zeta{(2\,n+1)}$, $\,n\,$ being a positive integer, which contains a series they could not reduce to a finite closed-form.  This open problem reveals a greater complexity in comparison to $\zeta(2n)$, which is a rational multiple of $\pi^{2n}$. For the Dirichlet beta function, the things are `inverse': $\beta(2n+1)$ is a rational multiple of $\pi^{2n+1}$ and no closed-form expression is known for $\beta(2n)$.  Here in this work, I modify the Dancs-He approach in order to derive an Euler-type formula for $\,\beta{(2n)}$, including $\,\beta{(2)} = G$, the Catalan's constant. I also convert the resulting series into zeta series, which yields new exact closed-form expressions for a class of zeta series involving $\,\beta{(2n)}$ and a finite number of odd zeta values. A closed-form expression for a certain zeta series is also conjectured.
\newline
\end{abstract}


\end{frontmatter}

\section{Introduction}

For real values of $s$, $\,s>1$, the Riemann zeta function is defined as $\zeta(s) := \sum_{n=1}^\infty{{\,1/n^s}}$. In this domain, the convergence of this series is guaranteed by the integral test.\footnote{For $s=1$, it represents the harmonic series $\sum_{n=1}^\infty{{\,1/n}}$, which diverges to infinity.}  For integer values of $s$, it is known that $\zeta{(2 n)}$, $n$ being a positive integer, the so-called even zeta values are rational multiples of $\,\pi^{2n}$, according to a famous Euler's formula (see, e.g.,~\cite{Tsumura} and references therein), namely
\begin{equation}
\zeta(2n) = (-1)^{n-1} \, \frac{2^{2n-1} \, B_{2n}}{(2n)!} \; \pi^{2n} \, ,
\label{eq:Euler}
\end{equation}
in which $n$ is a positive integer and $B_{n}$ are the Bernoulli numbers, i.e. the rational coefficients of ${\,z^n/n!}\,$ in the Taylor expansion of ${\,z/(e^z-1)}$, $|z| < 2\,\pi$. From the fact that $\pi$ is a transcendental number, as first proved by Lindemann (1882), it follows that all $\zeta(2n)$ are transcendental numbers.  For $\,\zeta{(2n+1)}$, on the other hand, no analogous closed-form expression and no transcendence proof is known. In fact, even irrationality proofs are unknown, except for the Ap\'{e}ry proof (1978) that $\,\zeta{(3)}\,$ is irrational~\cite{Apery,Porten}.  The belief that all these numbers are irrational increased much when Rivoal presented a proof that infinitely many odd zeta values are irrational~\cite{Rivoal}. This proof was soon followed by the `finite' result by Zudilin (2001) that at least one of the four numbers $\zeta(5), \zeta(7), \zeta(9),$ and $\zeta(11)$ is irrational~\cite{Zudilin}. 

This mysterious scenario has a `reverse' counterpart composed by the values of the Dirichlet beta function at integer points. This function is defined as\footnote{The function $\beta{(s)}$ can also be defined in terms of the Hurwitz zeta function $\zeta{(s,a)} := \sum_{k=0}^{\infty} 1/(k+a)^s$, as given by $\beta(s) = {\,\left[\zeta{(s,\frac14)}-\zeta{(s,\frac34)}\right]/4^s}$. \smallskip }
\begin{equation}
\beta(s) := \sum_{k=0}^\infty{\frac{(-1)^k}{(2\,k+1)^s}} \, ,
\label{eq:Bdef}
\end{equation}
a series that converges for all $s>0$, according to the Leibnitz's test for alternating series.  For this function, one has the following analogue of Eq.~\eqref{eq:Euler}:
\begin{equation}
\beta(2n+1) = (-1)^{n} \, \frac{E_{2n}}{2^{2n+2} \, (2n)!} \; \pi^{2n+1} \, ,
\label{eq:Euler2}
\end{equation}
where $E_{n}$ are the Euler numbers, i.e. the integer numbers obtained as the coefficients of ${\,z^n/n!}\,$ in the Taylor expansion of $\mathrm{sech}{(z)}$, $|z| < {\,\pi/2}\,$.\footnote{The analytic continuation of $\beta(s)$ extends it to all points in the complex plane, without singularities. This is done with the functional equation $\,\beta(1-s) = ({\,2/\pi})^{s-1} \, \sec{\left(\frac{\pi}{2} \, s\right)} {\: \beta(s) / \Gamma(s-1)}$, where $\Gamma{(x)}$ is the Euler's gamma function. \smallskip}  
Of course, the transcendence of $\pi$ implies the transcendence of all odd beta values, but no irrationality proof is known for any $\,\beta{(2 n)}$, not even for $\beta{(2)} = G = 0.91596559\ldots\,$, known as the Catalan's constant,\footnote{The suspected \emph{irrationality} of $G$ remains unproven. In fact, it is one of the main open problems in analytic number theory.}  which is the beta counterpart of $\zeta{(3)}$~\cite{bookCtes}.  Presently, the only known irrationality results are the recent proofs by Rivoal and Zudilin (2003) that infinitely many $\beta{(2n)}$ are irrational, and that at least one of the seven numbers $\beta(2), \, \ldots , \beta(14)$ is irrational~\cite{RZ}. Furthermore, none knows a finite closed-form expression for any $\beta(2n)$ in terms of elementary functions. However, in investigating some special values of the polygamma function $\psi^{(n)}{(x)}$ at rational entries,\footnote{Here, $\psi^{(n)}{(x)}$ is the $n$-th derivative of the digamma function $\psi(x)$, which in turn is defined as the logarithm derivative of $\Gamma{(x)}$. \smallskip}  K\"{o}lbig (1996) found a relation between $\psi^{(2n-1)}{({\,1/4})}$ and $\beta(2n)$~\cite{Kolbig}. By isolating $\beta(2n)$ in that relation, one has
\begin{equation}
\beta(2n) = \frac{\psi^{(2n-1)}{\left( \frac14 \right)}}{2\,(2n-1)!\,4^{2n-1}} \, - \frac{(2^{2n}-1) \, |B_{2n}|}{ 2\, (2n)!} \, \pi^{2n} \, ,
\label{eq:Kolbig}
\end{equation}
which resembles the Euler's formula. Unfortunately, the arithmetic nature of $\psi^{(2n-1)}{({\,1/4})}$ is also unknown and it seems very hard to express these numbers in terms of other known constants.

Apart from this singular result, the approach introduced by Dancs and He in a recent work~\cite{Dancs}, in which an Euler-type formula is derived for $\,\zeta{(2n+1)}\,$, suggests that a similar Euler-type formula could be found for $\beta(2n)$.  Here in this work, I make some series manipulations similar to those carried out by Dancs and He with a view to deriving an Euler-type formula for $\beta(2n)$. As will be shown, this yields a formula containing a series involving rational multiples of even powers of $\pi$, the companying coefficients involving the numbers $\,E_{\,2\,n+1}(1)$, where $\,E_{n}(x)$ is the Euler polynomial of degree $n$.\footnote{The Euler polynomials are $\:E_{n}(x) := \sum_{\,k=0}^{\,n} {\binom{n}{k} \dfrac{E_k}{2^k} \, \left( x -\frac12 \right)^{n-k}}$. \smallskip}  This formula can then be regarded as the analogue of Eq.~(5) of~\cite{Tsumura}, in the sense that the series reduces to a single term when $E_{2n+1}(1)$ is substituted by $E_{2n}(1)$. Lastly, by converting the summand into rational multiples of $\,\zeta{(2n)}$ and then making use of a formula derived recently by Milgran~\cite{Milgran}, I derive an exact closed-form expression for a certain class of zeta series related to $\beta{(2n)}$ and a finite number of odd zeta values.

\section{Adapting the Dancs-He approach for even beta values}

Let us modify the method introduced by Dancs and He, which yields an Euler-type formula for $\ln{2}$ and $\zeta{(2n+1)}$, $n$ being a positive integer,\footnote{In this context, $\,\ln{2}\,$ plays the role of $\zeta(1)$, the only singularity of $\zeta(s)$ in the complex plane. \smallskip } in order to get a similar result for $\beta{(2n)}$.  For a given $\epsilon > 0$ and $u \in [1,1+\epsilon]$, we start with the following series expansion
\begin{equation}
\frac{2 \, e^t}{e^t + u} = \sum_{n=0}^{\infty}{\phi_n(u)\, \frac{t^n}{n!}} \,.
\label{eq:phiserie}
\end{equation}
From the generating function for $E_n(x)$, namely ${\,2 \, e^{x\,t}/(e^t + 1)} = \sum_{n=0}^{\infty}{E_n(x)\, \frac{t^n}{n!}}$, it is clear that $\phi_n(1) = E_n(1)$, $\forall \, n \ge 0$. For $u>1$, we have
\begin{equation}
\phi_n(u) = -2 \, \sum_{j=1}^{\infty}{\frac{j^{\,n}}{(-u)^{\,j}}} \, .
\end{equation}
For a given $u \ge 1$, this series converges absolutely for $|t| < \pi\,$~\cite{Dancs}. For $n<0$, we shall take this series as our definition of $\phi_n(u)$. From the Taylor series for $\ln{(1+x)}$, we have $\phi_{-1}(1) = 2\,\ln{2}$, and from a well-known result for the alternating zeta function, we have
\begin{equation}
\phi_{-m}(1) = 2\,(1-2^{1-m})\,\zeta{(m)}
\label{eq:phi}
\end{equation}
for all $\: m \in \mathbb{Z}$, $\,m>1$.

Now, let us modify the argument of the sine function in the first series that appear after Eq.~(2.4) in~\cite{Dancs} by exchanging $\,\pi\,$ by ${\,\pi/2\,}$. This results in
\begin{eqnarray*}
\sum_{n=1}^\infty{ \frac{1}{(-u)^n} \, \frac{\sin{\left(n\,\frac{\pi}{2}\right)}}{n^{2k}} } &=& \sum_{n=1}^\infty{ (-1)^n \, \frac{\sin{\left(n\,\frac{\pi}{2}\right)}}{u^n \, n^{2k}} }  \\
&=& \sum_{n \: \mathrm{even}}{ (-1)^n \, \frac{\sin{\left(n\,\frac{\pi}{2}\right)}}{u^n \, n^{2k}} } +\sum_{n \: \mathrm{odd}}{ (-1)^n \, \frac{\sin{\left(n\,\frac{\pi}{2}\right)}}{u^n \, n^{2k}} } \\
&=& \sum_{m=1}^\infty{ (-1)^{2m} \, \frac{\sin{(m\,\pi)}}{u^{2m} \, (2m)^{2k}} } +\sum_{m=1}^\infty{ (-1)^{2m-1} \, \frac{\sin{\left(m\,\pi -\frac{\pi}{2}\right)}}{u^{2m-1} \, (2m-1)^{2k}} } \,.
\end{eqnarray*}
The series with $\,\sin{(m \pi)}\,$ is of course null and, by noting that $\,\sin{\!\left(m\,\pi -\frac{\pi}{2}\right)} = (-1)^{m+1}$, one has
\begin{equation}
\sum_{n=1}^\infty{ (-1)^n \, \frac{\sin{\left(n\,\frac{\pi}{2}\right)}}{u^n \, n^{2k}} }  = \sum_{m=1}^\infty{ (-1)^m \, \frac{1}{u^{2m-1} \, (2m-1)^{2k}} } \, .
\label{eq:mestra}
\end{equation}
By taking the limit as $u \rightarrow 1^{^{+}}$ on both sides of this equation, one finds, for any positive integer $k$,
\begin{equation}
\lim_{u \rightarrow 1^{^{+}}}{\:\sum_{n=1}^\infty{ (-1)^n \, \frac{\sin{\left(n\,\frac{\pi}{2}\right)}}{u^n \, n^{2k}} }}  = \sum_{m=1}^\infty{\frac{(-1)^m}{(2m-1)^{2k}} } \, .
\label{eq:aux}
\end{equation}
From Eq.~\eqref{eq:Bdef}, one has
\begin{equation}
\lim_{u \rightarrow 1^{^{+}}}{\:\sum_{n=1}^\infty{ (-1)^n \, \frac{\sin{\left(n\,\frac{\pi}{2}\right)}}{u^n \, n^{2k}} }} = -\beta{(2k)} \, .
\label{eq:b2k}
\end{equation}
This limit can be calculated by manipulating the series in a manner to express it in terms of the functions $\phi_n(u)$. This is done by using the Taylor expansion of the sine function in the summand, as follows:
\begin{eqnarray*}
-\beta{(2k)} &=& \lim_{u \rightarrow 1^{^{+}}}{\:\sum_{n=1}^\infty{ \left[(-1)^n \, \frac{1}{u^n \, n^{2k}} \sum_{j=0}^\infty{(-1)^j \, \frac{\left(n\,\frac{\pi}{2}\right)^{2j+1}}{(2j+1)!} } \right] }} \\
&=& \lim_{u \rightarrow 1^{^{+}}}{\, \frac{u}{2} \sum_{j=0}^\infty{(-1)^{j+1} \, \frac{\pi^{2j+1}}{2^{2j+1}\,(2j+1)!} \: \phi_{2j+1-2k}(u) } } \\
&=& \lim_{u \rightarrow 1^{^{+}}}{\, \frac{u}{2} \left[ \, \sum_{j=0}^{k-1}{\frac{(-1)^{j+1} \, \pi^{2j+1}}{2^{2j+1}\,(2j+1)!} \: \phi_{2j+1-2k}(u) } + \sum_{j=k}^{\infty}{\frac{(-1)^{j+1} \, \pi^{2j+1}}{2^{2j+1}\,(2j+1)!} \: \phi_{2j+1-2k}(u) } \right]} \\
&=& \lim_{u \rightarrow 1^{^{+}}} \frac{u}{2} \, \Bigg[ \, \sum_{j=0}^{k-1}{\frac{(-1)^{j+1} \: ({\,\pi/2})^{2j+1}}{(2j+1)!} \: \phi_{2j+1-2k}(u) } \\
&+& \sum_{m=0}^\infty{\frac{(-1)^{m+k+1} \: ({\,\pi/2})^{2m+2k+1}}{(2m+2k+1)!} \: \phi_{2m+1}(u) } \Bigg] ,
\end{eqnarray*}
where the last series was obtained by substituting $m = j-k$. Lastly, by effectively taking the limit as $u \rightarrow 1^{^{+}}$, one finds
\begin{eqnarray}
\beta{(2k)} &=& \frac12 \Bigg[ \, \sum_{j=0}^{k-1}{\frac{(-1)^j \: (\pi/2)^{2j+1}}{(2j+1)!} \: \phi_{2j+1-2k}(1) } \nonumber \\
&+& (-1)^k \left(\frac{\pi}{2} \right)^{2k+1} \sum_{m=0}^\infty{\frac{f_m}{(2m+2k+1)! \cdot 2^{\,2m}} } \Bigg] ,
\label{eq:geratriz}
\end{eqnarray}
where $f_m := (-1)^m \: {\pi}^{2m} \, E_{2m+1}(1)\,$ is the same parameter defined in Eq.~(2.5) of~\cite{Dancs}. With this equation in hands, we can now state and prove the following theorem.

\begin{teo}[Euler-type formula for $\beta{(2 n)}\,$]
 \label{teo:b2n}
\; ~For any integer $k$, $k>0$,
\begin{eqnarray}
\beta{(2k)} &=& (-1)^{k+1}\,\frac{({\,\pi/2})^{2k-1}}{(2k-1)!} \, \ln{2} \, - (-1)^{k} \, \sum_{m=1}^{k-1} \frac{(-1)^m \: ({\,\pi/2})^{2k-2m-1}}{(2k-2m-1)!} \left(1-\frac{1}{2^{2m}}\right) \nonumber \\
&\times& \zeta(2m+1) +(-1)^k \, \frac{\pi^{2k+1}}{2^{2k+2}} \, \sum_{m=0}^\infty{(-1)^m \frac{{\pi}^{2m} \: E_{2m+1}(1)}{(2m+2k+1)! \cdot 2^{\,2m}} } \: .
\label{eq:teo1}
\end{eqnarray}
\end{teo}

\,

\begin{prova}
\quad By making explicit the last term in the finite sum in Eq.~\eqref{eq:geratriz}, one has
\begin{eqnarray*}
\beta{(2k)} &=& \frac12 \, \Bigg[ \sum_{j=0}^{k-2}{\frac{(-1)^j \: (\pi/2)^{2j+1}}{(2j+1)!} \: \phi_{2j-2k+1}(1) } \, + (-1)^{k+1} \, \frac{({\,\pi/2})^{2k-1}}{(2k-1)!} \: \phi_{-1}(1) \\
&+& (-1)^k \, \left(\frac{\pi}{2} \right)^{2k+1} \, \sum_{m=0}^\infty{\frac{f_m}{(2m+2k+1)! \cdot 2^{\,2m}} } \Bigg] .
\end{eqnarray*}
By substituting both $\,\phi_{-1}(1) = 2\,\ln{2}\,$ and the result for $\phi_{-m}(1)$ stated in Eq.~\eqref{eq:phi}, one has
\begin{eqnarray*}
\beta{(2k)} &=& \frac12 \, \Bigg[ 2 \sum_{j=0}^{k-2}{\frac{(-1)^j \: (\pi/2)^{2j+1}}{(2j+1)!} \: \left(1-2^{2j-2k+2}\right) \cdot \zeta(2k-2j-1)} \\
&+& 2\,(-1)^{k+1} \, \frac{({\,\pi/2})^{2k-1}}{(2k-1)!} \, \ln{2} \: + (-1)^k \, \left(\frac{\pi}{2} \right)^{2k+1} \, \sum_{m=0}^\infty{\frac{f_m}{(2m+2k+1)! \cdot 2^{\,2m}} } \Bigg] .
\end{eqnarray*}
This simplifies to
\begin{eqnarray*}
\beta{(2k)} &=& \sum_{j=0}^{k-2}{\frac{(-1)^j \: (\pi/2)^{2j+1}}{(2j+1)!} \: \left(1-2^{2j-2k+2}\right) \cdot \zeta(2k-2j-1)} \, \\
&+& (-1)^{k+1} \, \frac{({\,\pi/2})^{2k-1}}{(2k-1)!} \, \ln{2} + (-1)^k \, \frac{\pi^{2k+1}}{2^{2k+2}} \, \sum_{m=0}^\infty{\frac{f_m}{(2m+2k+1)! \cdot 2^{\,2m}} } \, .
\end{eqnarray*}
The substitution $j = k-m-1$ in the finite sum completes our proof.
\begin{flushright} $\Box$ \end{flushright}
\end{prova}

This new formula for $\beta{(2 k)}$ can be regarded as the analogue of Eq.~(5) of~\cite{Tsumura}. 
This may give some insight into why it is so difficult to find closed-form expressions for even beta values.

Let us list the result for the first three beta values yielded by Theorem~\ref{teo:b2n}, as found by putting $k=1, 2, 3$ on Eq.~\eqref{eq:teo1}:
\begin{equation}
\beta{(2)} = G = \frac{\pi}{2} \, \ln{2} - \frac{\pi^3}{2^4} \, \sum_{m=0}^\infty{(-1)^m \frac{{\pi}^{2m} \: E_{2m+1}(1)}{(2m+3)! \cdot 2^{\,2m}} } \, ,
\label{eq:G}
\end{equation}

\begin{equation}
\beta{(4)} = -\,\frac{\pi^3}{3! \cdot 2^3} \, \ln{2} + \frac{3}{2^3} \, \pi \, \zeta{(3)} + \frac{\pi^5}{2^6} \, \sum_{m=0}^\infty{(-1)^m \frac{{\pi}^{2m} \: E_{2m+1}(1)}{(2m+5)! \cdot 2^{\,2m}} } \, ,
\label{eq:B4}
\end{equation}
and
\begin{equation}
\beta{(6)} = \frac{\pi^5}{5! \cdot 2^5} \, \ln{2} - \frac{\pi^3}{2^6} \, \zeta{(3)} + \frac{15\,\pi}{2^5} \, \zeta{(5)} - \frac{\pi^7}{2^8} \, \sum_{m=0}^\infty{(-1)^m \frac{{\pi}^{2m} \: E_{2m+1}(1)}{(2m+7)! \cdot 2^{\,2m}} } \, ,
\label{eq:B6}
\end{equation}
respectively.

\section{A closed-form expression for a class of zeta series}

Let us express the infinite series in the previous formulas in the form of zeta series, in order to establish a finite closed-form expression for a class of zeta series related to $\beta{(2k)}$~\cite{Sirivas}. For this, let us make use of the identities
\begin{equation}
E_{2m+1}(1) = -E_{2m+1}(0) = 2 \, \frac{2^{2m+2}-1}{2m+2}\,B_{2m+2} \, .
\end{equation}
By making use of the Euler's formula, see Eq.~\eqref{eq:Euler}, one finds the following practical conversion formula:
\begin{equation}
f_m = (-1)^m \, \pi^{2m} \, E_{2m+1}(1) = \frac{4-2^{-2m}}{\pi^2} \, (2m+1)! \cdot \zeta{(2m+2)} \, .
\label{eq:converte}
\end{equation}

By substituting this on Eq.~\eqref{eq:G} and then changing the summation index $m$ by $n-1$, one finds
\begin{equation}
G = \frac{\pi}{2} \, \ln{2} - \pi \left[\: \sum_{n=1}^\infty{\frac{\zeta(2n)}{2n\,(2n+1)\, 2^{2n}}} \, -\sum_{n=1}^\infty{\frac{\zeta(2n)}{2n\,(2n+1)\, 4^{2n}}} \right] .
\label{eq:Gzeta}
\end{equation}
The first series evaluates to ${\,(\ln{\pi}-1)/2}$, accordingly to a recent result by Fujii and Suzuki (see Eqs.(9--11) in~\cite{Fuji}). The remaining zeta series then simplifies to
\begin{equation}
\sum_{n=1}^\infty{\frac{\zeta(2n)}{n\,(2n+1)\, 4^{2n}}} = 2\, \frac{G}{\pi} + \ln{\left(\frac{\pi}{2}\right)}-1 \, ,
\label{eq:newGzeta}
\end{equation}
a result that can be found on p.~242, Eq.~(672) of~\cite{Sirivas}. It can also be deduced from a series involving the Clausen function $\mathrm{Cl}_2(\theta)$, at p.~265 of~\cite{revBorwein}, by using the special value $\mathrm{Cl}_2(\pi/2) = G$.

As another example, let us apply our conversion formula, Eq.~\eqref{eq:converte}, in the expression for $\beta{(4)}$ established in Eq.~\eqref{eq:B4}. One finds
\begin{eqnarray}
\beta{(4)} &=& -\,\frac{\pi^3}{48} \, \ln{2} + \frac38 \, \pi \, \zeta{(3)} + \frac{\pi^3}{4} \, \Bigg[ \, \sum_{n=1}^\infty{\frac{\zeta(2n)}{2n \, (2n+1)\,(2n+2)\,(2n+3)\, 2^{2n}}} \nonumber \\
&-&\sum_{n=1}^\infty{\frac{\zeta(2n)}{2n \, (2n+1)\,(2n+2)\,(2n+3)\, 4^{2n}}} \, \Bigg] .
\label{eq:B4zeta}
\end{eqnarray}
The first series evaluates to ${\,\zeta(3)/(2 \pi^2)} +{\,\ln{\pi}/12} -{\,11/72}$, accordingly to Wilton's formula (see Eq.~(31) at p.~148 of~\cite{Sirivas}, also Eq.~(54) in~\cite{Cho2006}), which yields
\begin{equation}
\sum_{n=1}^\infty{\frac{\zeta(2n)}{2n \, (2n+1)\,(2n+2)\,(2n+3)\, 4^{2n}}} = -4 \, \frac{\beta{(4)}}{\pi^3} + 2 \,\frac{\zeta(3)}{\pi^2} +\frac{1}{12} \, \ln{\left(\frac{\pi}{2} \right)} -\frac{11}{72} \, ,
\label{eq:newB4zeta}
\end{equation}
This series is equivalent to that obtained by putting $\,t = \frac14\,$ in
\begin{eqnarray}
\sum_{k=1}^\infty{\frac{\zeta(2k)}{k \, (k+1)\,(2k+1)\,(2k+3)} \, t^{2k+3}} = \frac{\zeta(3)}{2\,\pi^2}\, t + \frac{6\,\ln{2\pi} -11}{18}\,t^3 \nonumber \\
+\frac13 \left[ \zeta'(-3,1+t) -\zeta'(-3,1-t) \right] , \quad |t| < 1 \, ,
\end{eqnarray}
which is one of the zeta series investigated by Srivastava and co-workers (see Eq.~(713) in~\cite{siri}). Here, $\zeta'(s,a)$ denotes the first derivative with respect to $s$.

As a last example, let us derive the zeta series corresponding to $\beta{(6)}$. By applying our conversion formula to the expression for $\beta{(6)}$ established in Eq.~\eqref{eq:B6}, one finds, after some algebra, that
\begin{eqnarray}
\sum_{n=1}^\infty{\frac{\zeta(2n)}{2n \, (2n+1)\,(2n+2)\,(2n+3)\,(2n+4)\,(2n+5)\,} \cdot \left(\frac{1}{2^{2n}}-\frac{1}{4^{2n}}\right)} = \nonumber \\
\frac{1}{240}\,\ln{2} -\frac14 \frac{\zeta(3)}{\pi^2} +\frac{15}{2} \frac{\zeta(5)}{\pi^4} -16 \, \frac{\beta{(6)}}{\pi^5} \, .
\label{eq:B6zetaAux}
\end{eqnarray}
By making use of the software Mathematica (release 7), one finds the following exact analytical result:
\begin{eqnarray}
\sum_{n=1}^\infty{\frac{\zeta(2n)}{2n \, (2n+1)\,(2n+2)\,(2n+3)\,(2n+4)\,(2n+5)\,} \cdot \frac{1}{2^{2n}}} = \nonumber \\
\frac{\zeta(3)}{12\,\pi^2} -\frac{\zeta(5)}{2\,\pi^4} +\frac{\ln{\pi}}{240} -\frac{137}{14400}  \, .
\label{eq:B6zetaAux2}
\end{eqnarray}
By substituting this on Eq.~\eqref{eq:B6zetaAux}, one finds
\begin{eqnarray}
\sum_{n=1}^\infty{\frac{\zeta(2n)}{2n \, (2n+1)\,(2n+2)\,(2n+3)\,(2n+4)\,(2n+5)\,} \cdot \frac{1}{4^{2n}}} = \nonumber \\
\frac{\zeta(3)}{3\,\pi^2} -8\,\frac{\zeta(5)}{\pi^4} +16\,\frac{\beta{(6)}}{\pi^5} +\frac{\ln{\left(\frac{\pi}{2}\right)}}{240} -\frac{137}{14400} \, ,
\label{eq:B6zetaAux3}
\end{eqnarray}
a closed-form result that is not found in literature.


Let us generalize the above results for zeta series involving $\beta(2\,k)$.

\begin{teo}[Zeta series involving $\beta{(2k)}\,$]
 \label{teo:zetab2n}
\quad For any positive integer $k$,
\begin{eqnarray}
\sum_{n=1}^\infty{\frac{\zeta(2n)}{2n \, (2n+1) \, \ldots \, (2n+2k-1)} \cdot \left(\frac{1}{2^{2n}}-\frac{1}{4^{2n}}\right)} = \nonumber \\
(-1)^k\,\frac{2^{2k-2}}{\pi^{2k-1}} \, \beta{(2k)} 
+\frac{k}{(2k)!}\,\ln{2} + \, \frac12 \, \sum_{m=1}^{k-1}{(-1)^m \, \frac{2^{2m}-1}{\pi^{2m}} \,  \frac{\zeta{(2m+1)}}{(2k-2m-1)!}} \, .
\label{eq:teo2}
\end{eqnarray}
\end{teo}

\begin{prova}
\quad  Take the Euler-type formula for $\beta{(2k)}$ established in Theorem~\ref{teo:b2n}. By using our conversion formula, Eq.~\eqref{eq:converte}, one finds, after simplifications,
\begin{eqnarray*}
\beta{(2k)} &=& (-1)^{k+1}\,\frac{({\,\pi/2})^{2k-1}}{(2k-1)!} \, \ln{2} \\
&+& (-1)^{k+1} \, \sum_{m=1}^{k-1} (-1)^m \, \frac{({\,\pi/2})^{2k-2m-1}}{(2k-2m-1)!} \left(1-\frac{1}{2^{2m}}\right) \cdot \zeta(2m+1) \\
&+& (-1)^k \, \frac{\pi^{2k-1}}{2^{2k}} \, \sum_{m=0}^\infty{\frac{\zeta(2m+2)}{(2m+2)\,(2m+3)\,\ldots\,(2m+2k+1)} \, \left(\frac{1}{2^{2m}}-\frac{1}{2^{4m+2}}\right)} .
\end{eqnarray*}
By doing $m=n-1$ in the infinite series, one has
\begin{eqnarray*}
\beta{(2k)} &=& (-1)^{k+1} \frac{({\pi/2})^{2k-1}}{(2k-1)!} \, \ln{2} \\
&-& (-1)^{k} \sum_{m=1}^{k-1} (-1)^m \frac{({\,\pi/2})^{2k-2m-1}}{(2k-2m-1)!} \left(1-\frac{1}{2^{2m}}\right) \cdot \zeta(2m+1) \\
&+& (-1)^k \, 2\, \left(\frac{\pi}{2}\right)^{2k-1} \, \sum_{n=1}^\infty{\frac{\zeta(2n)}{2n\,(2n+1)\,\ldots\,(2n+2k-1)} \, \left(\frac{1}{2^{2n}}-\frac{1}{2^{4n}}\right)} .
\end{eqnarray*}
A multiplication by $(-1)^k \, ({\,2/\pi})^{2k-1}$ on both sides, followed by the isolation of the zeta series, completes the proof.
\begin{flushright} $\Box$ \end{flushright}
\end{prova}

As expected, our previous closed-form expressions for zeta series stated in Eqs.~\eqref{eq:Gzeta},~\eqref{eq:B4zeta}, and~\eqref{eq:B6zetaAux}, are special cases ($k=1,2,3$,  respectively) of the general formula in Eq.~\eqref{eq:teo2}. Interestingly, this new closed-form expression for the general zeta series can be regarded as the $\beta(2k)$ counterpart of a formula derived by Milgran (2007) for a zeta series related to $\zeta{(2k+1)}$~\cite{Milgran}.

Clearly, our zeta series involving $\beta{(2k)}$, as found in Eq.~\eqref{eq:teo2}, can be written as the difference of two zeta series, as done in the derivation of Eqs.~\eqref{eq:newGzeta}, \eqref{eq:newB4zeta}, and \eqref{eq:B6zetaAux3}. This can be justified by noting that both the series
\begin{equation}
\sum_{n=1}^\infty{\frac{\zeta(2n)}{2n \, (2n+1) \, \ldots \, (2n+2k-1)} \cdot \frac{1}{2^{\,2n}}}
\label{eq:ser1}
\end{equation}
and
\begin{equation}
\sum_{n=1}^\infty{\frac{\zeta(2n)}{2n \, (2n+1) \, \ldots \, (2n+2k-1)} \cdot \frac{1}{4^{\,2n}}}\label{eq:ser2}
\end{equation}
converge absolutely. The task of finding a closed-form expression for one of these series demands the previous knowledge of a closed-form result for the other. However, I could not find a closed-form result for any of these series, neither in literature nor by using mathematical softwares (Mathematica and Maple). The best I could do was to investigate the pattern of the first few (exact) analytical results found for the zeta series in Eq.~\eqref{eq:ser1}.  This experimental procedure strongly indicates that
\begin{eqnarray}
\sum_{n=1}^\infty{\frac{\zeta(2n)}{2n \, (2n+1) \, \ldots \, (2n+N)} \cdot \frac{1}{2^{2n}}} = \nonumber \\
\frac12 \left[ \frac{\ln{\pi}}{N!} \, - \frac{H_N}{N!} + \sum_{m=1}^{(N-1)/2}{(-1)^{m+1} \, \frac{\zeta{(2m+1)}}{\pi^{2m}\,(N-2m)!}} \right] ,
\label{eq:conjec1}
\end{eqnarray}
where $H_N := \sum_{i=1}^N{1/i}$ is the $N$-th harmonic number, $N = 2 k-1$ being a positive odd integer. Although I could not find any exception to the validity of this analytical closed-form expression, at the current stage of my research on zeta series this is only a \emph{conjecture}.\footnote{I have also tested the validity of Eq.~\eqref{eq:conjec1} numerically for every positive integer value of $k$ up to $1000$, which corresponds to $N=1999$.}  Note that a formal proof for this conjecture will automatically imply, from our Theorem~2, that
\begin{eqnarray}
\sum_{n=1}^\infty{\frac{\zeta(2n)}{2n \, (2n+1) \, \ldots \, (2n+N)} \cdot \frac{1}{4^{2n}}} = \nonumber \\
\frac12 \left[ \frac{\ln{(\pi/2)}}{N!} -\frac{H_N}{N!} -(-1)^{\frac{N+1}{2}} \left(\frac{2}{\pi}\right)^N \beta{(N+1)} -\sum_{m=1}^{\frac{N-1}{2}}{(-1)^m \left(\frac{2}{\pi}\right)^{2m} \frac{\zeta{(2m+1)}}{(N-2m)!}} \right]
\label{eq:conjec2}
\end{eqnarray}
is also true.

\section*{Acknowledgments}
Thanks are due to Mrs. Marcia R. Souza for checking the convergence and validity of all series expansions in this work.


\begin{thebibliography}{99}

\bibitem{Apery} R. Ap\'{e}ry, \emph{Irrationalit\'{e} de $\zeta(2)$ et $\zeta(3)$}, Ast\'{e}risque \textbf{61} (1979), pp.~11--13.

\bibitem{revBorwein} J. M. Borwein, D. M. Bradley, and R. E. Crandall, \emph{Computational strategies for the Riemann zeta function}, J. Comput. Appl. Math. \textbf{121} (2000), pp.~ 247--296.

\bibitem{Cho2006} Y. J. Cho, M. Jung, J. Choi, and H. M. Srivastava, \emph{Closed-form evaluations of definite integrals and associated infinite series involving the Riemann zeta function}, Int. J. Comp. Math. \textbf{83} (2006), pp.~461--472.

\bibitem{Dancs} M. J. Dancs and T.-X. He, \emph{An Euler-type formula for $\zeta{(2k+1)}$}, J. Number Theory \textbf{118} (2006), pp.~192--199.

\bibitem{bookCtes} S. R. Finch, \emph{Mathematical Constants}, Cambridge Univ. Press, Cambridge UK, 2003.

\bibitem{Fuji} K. Fujii and T. Suzuki, \emph{Introduction of an elementary method to express $\zeta(2n+1)$ in terms of $\zeta(2k)$ with $k \ge 1$}. \, Available at arXiv: 0805.0030v2 (2008).

\bibitem{Kolbig} K. S. K\"{o}lbig, \emph{The polygamma function $\psi^{(k)}(x)$ for $x =1/4$ and $x =3/4$}, J. Comput. Appl. Math. \textbf{75} (1996), pp.~43--46.

\bibitem{Milgran} M. Milgran, \emph{Notes on a paper of Tyagi and Holm: ``A new integral representation for the Riemann zeta function''}. \, Available at  arXiv: 0710.0037v1  (2007).

\bibitem{Porten} A. van der Poorten, \emph{A proof that Euler missed: Ap\'{e}ry's proof of the irrationality of $\,\zeta(3)$}, Math. Intelligencer \textbf{1} (1979), pp.~195--203.

\bibitem{Rivoal} T. Rivoal, \emph{La fonction z\^{e}ta de Riemann prend une infinit\'{e} de valeurs irrationnelles aux entiers impairs}, C. R. Acad. Sci. Paris \textbf{331}, S\'{e}rie I, (2000), pp.~267--270.

\bibitem{RZ} T. Rivoal and W. Zudilin, \emph{Diophantine properties of numbers related to Catalan's constant}, Math. Annalen \textbf{326} (2003), pp.~705--721.

\bibitem{siri}  H. M. Srivastava, M. L. Glasser, and V. S. Adamchik, \emph{Some definite integrals associated with the Riemann zeta function}, Zeitschrift f\"{u}r Analysis and Anwendungen \textbf{19} (2000), pp.~831--846.

\bibitem{Sirivas}  H. M. Srivastava and J. Choi, \emph{Series associated with the zeta and related functions}, Kluwer Academic Publishers, Dordrecht, 2001.

\bibitem{Tsumura} H. Tsumura, \emph{An elementary proof of Euler's formula for $\zeta{(2\,m)}$}, Amer. Math. Monthly \textbf{111} (2004), pp.~430--431.

\bibitem{Zudilin} W. Zudilin, \emph{One of the numbers $\zeta(5)$, $\zeta(7)$, $\zeta(9)$, $\zeta(11)$ is irrational}, Russian Math. Surveys \textbf{56} (2001), pp.~774--776.

\end{thebibliography}
\end{document}